\documentclass[12pt, twoside, leqno]{article}
\usepackage{amsmath,amsthm}
\usepackage{amssymb}
\usepackage{enumerate}
\usepackage{graphicx}

\newtheorem{thm}{Theorem}[section]

\theoremstyle{definition}

\def\Q{{\mathbb Q}}
\def\Z{{\mathbb Z}}

\def\O{{\cal O}}

\numberwithin{equation}{section}

\begin{document}
\baselineskip=17pt

\title{
Integral bases and monogenity of the \\simplest sextic fields\\
}
\author{
Istv\'{a}n Ga\'{a}l\thanks{
        Research supported in part by K115479 from the
        Hungarian National Foundation for Scientific Research				
				and by the EFOP-3.6.1-16-2016-00022 project. 
				The project is co-financed by the European Union and the European Social Fund.	
                         },\; 
and L\'aszl\'o Remete
\thanks{
        Research supported in part through the new 
				national excellence program of the ministry of human capacities.
                         }
\\ \\
University of Debrecen, Mathematical Institute \\
H--4002 Debrecen Pf.400., Hungary \\
e--mail: gaal.istvan@unideb.hu, remete.laszlo@science.unideb.hu
}
\date{}

\maketitle

\renewcommand{\thefootnote}{}



\renewcommand{\thefootnote}{\arabic{footnote}}
\setcounter{footnote}{0}

\begin{abstract}
Let $m$ be an integer, $m\neq -8,-3,0,5$ such that $m^2+3m+9$ is square free.
Let $\alpha$ be a root of
\[
f=x^6-2mx^5-(5m+15)x^4-20x^3+5mx^2+(2m+6)x+1.
\]
The totally real cyclic fields $K=\Q(\alpha)$ are called simplest sextic fields
and are well known in the literature.

Using a completely new approach we explicitly give an integral basis 
of $K$ in a parametric form and 
we show that the structure of this integral basis is periodic in $m$ with period length 36.
We prove that $K$ is not monogenic except for a few values of $m$
in which cases we give all generators of power integral bases.
\end{abstract}

\section{Introduction}
\label{iiii}

Let $K$ be an algebraic number field of degree $n$ with ring of integers $\Z_K$.
It is a classical problem of algebraic number theory 
(cf. W.Narkiewicz \cite{nark}, I.Ga\'al \cite{gaal})
to decide 
if $\Z_K$ is {\bf monogenic}, that is if there exists an algebraic integer
$\alpha$ such that $\Z_K=\Z[\alpha]$. In this case 
 $\{1,\alpha,\ldots,\alpha^{n-1}\}$
is an integral basis called {\bf power integral basis}.
The coefficients of such an $\alpha$ in an integral basis of $K$
satisfies the so called {\bf index form equation} (see below).
Therefore in order to determine 
generators of power integral bases
authors were trying to develop algorithms for
solving index form equations.

There are efficient algorithms for solving index form equations
 in lower degree number fields (see \cite{gaal}),
mainly for cubic and quartic number fields. Some tedious algorithms
are known to solve the problem in quintic, sextic and some special octic, nonic 
fields. 

It is a much harder task to consider the problem of monogenity 
in {\bf infinite parametric families of number fields}. Let $m$ be an integer
parameter, let $f_m(x)\in\Z[x]$ be irreducible polynomials of degree $n$,
the coefficients of which depend on $m$.
Let $\alpha_m$ be a root of $f_m(x)$. We ask if $K_m=\Q(\alpha_m)$
is monogenic or not, moreover, what are the generators of power integral
bases of $K_m$?

It was possible to solve this parametric problem in some lower degree familes of
number fields, see e.g. \cite{gp}, \cite{gaal}. 
For parametric families of higher
degree number fields one of the main problems was 
that no integral basis of $K_m$ was known in a parametric form.
Therefore the problem of monogenity was often
considered in the {\bf equation
order} $\Z[\alpha_m]$, see e.g. \cite{gop}, instead of the ring of integers $\Z_K$.
Moreover, in higher degree number fields it is hopeless to solve the
index form equation in a parametric form because of the high degree 
($n(n-1)/2$) and the number of variables ($n-1$).

On the other hand our experience shows that higher degree number fields
seldom admit power integral bases. It turned out that in some cases it is 
much easier to show that $K_m$ is {\bf not monogenic} then to solve the
index form equation and to find that there are no solutions.

Recently the authors developed \cite{gr2} a completely new approach
to deal with {\bf infinite parametric families of higher degree
number fields}. We succeeded to determine an {\bf integral basis in a parametric
form}, managed to explicitly determine the factors of the corresponding index form
and proved that some linear combinations of these
factors are divisible by certain constants. These divisibility
conditions usually do not hold, therefore the fields are not monogenic.

Let $n$ be fixed with $2\leq n\leq 9$.
In \cite{gr2} we applied this method to pure fields 
with 
\[
f_m(x)=x^n-m,
\] 
where $m$ is square-free. 
It also turned out that the structure 
of an integral basis of $K_m=\Q(\sqrt[n]{m})$ is {\bf periodic} in $m$ which was
a new phenomenon.   

In this paper we develop further the method of 
\cite{gr2} and would like to demonstrate that it 
is applicable also to some other infinite parametric
families of number fields. 
In this paper we consider the family
\[
f_m(x)=x^6-2mx^5-(5m+15)x^4-20x^3+5mx^2+(2m+6)x+1,
\]
where $m$ is an integer $m\not=-8,-3,0,5$. If $\alpha_m$ is
a root of $f_m(x)$ then the fields $K_m=\Q(\alpha_m)$ are called 
{\bf simplest sextic fields}.  

\vspace{0.5cm}

\noindent
In the present paper \\
--we explicitly give an integral basis 
of $K_m$, \\
--we show that its structure is periodic with period length 36,\\
--we prove that $K_m$ is not monogenic except for $m=-4,-2,-1,1$,\\
--for $m=-4,-2,-1,1$ we present all generators of power integral bases.

\vspace{0.5cm}

We emphasize that to achieve these results we 
extended our method of \cite{gr2} at several points.
One of the main inventions is the application of the {\bf Dedekind
test} to construct an integral basis.
We hope that this way our method will find further applications, as well.

\section{Power integral bases and \\ monogenity of number fields}

Here we shortly recall those concepts (cf. \cite{gaal} for details)
that we use throughout.
Let $\alpha$ be a primitive integral element of the number field $K$
(that is $K=\Q(\alpha)$) of degree $n$ with ring of integers $\Z_K$. 
The {\bf index} of $\alpha$ is
\[
I(\alpha)=(\Z_K^{+}:\Z[\alpha]^{+})=\sqrt{\left|\frac{D(\alpha)}{D_K}\right|}
=\frac{1}{\sqrt{|D_K|}}\prod_{1\leq i<j\leq n}\left|\alpha^{(i)}-\alpha^{(j)}\right| \;\; ,
\]
where $D_K$ is the discriminant of $K$ and $\alpha^{(i)}$ denote
the conjugates of $\alpha$. 

If $B=\{b_1=1,b_2,\ldots,b_n\}$ is an
integral basis of $K$, then the {\bf index form}
corresponding to this integral basis is defined by
\[
I(x_2,\ldots,x_n)=\frac{1}{\sqrt{|D_K|}}
\prod_{1\leq i<j\leq n}
\left(
(b_2^{(i)}-b_2^{(j)})x_2+\ldots+(b_n^{(i)}-b_n^{(j)})x_n
\right)
\]
(where $b_j^{(i)}$ denote the conjugates of $b_j$).
This is a homogeneous polynomial with integral coefficients
of degree $n(n-1)/2$ in $n-1$ variables.
For any integral element 
\[
\alpha=x_1+b_2x_2+\ldots+b_nx_n
\]
(with $x_1,\ldots,x_n\in\Z$) we have
\[
I(\alpha)=|I(x_2,\ldots,x_n)|
\]
independently of $x_1$.
The element 
$\alpha$ generates a {\bf power integral basis} $\{1,\alpha,\ldots,\alpha^{n-1}\}$
if and only if $I(\alpha)=1$ that is if $(x_2,\ldots,x_n)\in \Z^{n-1}$ is 
a solution of the
{\bf index form equation}
\begin{equation}
I(x_2,\ldots,x_n)=\pm 1 \;\;\; {\rm in} \;\;\; (x_2,\ldots,x_n)\in \Z^{n-1}.
\label{iiixxx}
\end{equation}
In this case the ring of integers of $K$ is a simple ring extension
of $\Z$, that is $\Z_K=\Z[\alpha]$ 
and $K$ is called {\bf monogenic}.

\section{Simplest sextic fields}

Let $m$ be an integer $m\neq -8,-3,0,5$, and let $\alpha_m$ be a root of
\begin{equation}
f_m(x)=x^6-2mx^5-(5m+15)x^4-20x^3+5mx^2+(2m+6)x+1.
\label{f}
\end{equation}
The fields $K_m=\Q(\alpha_m)$ are called
{\bf simplest sextic fields} \cite{c6}, \cite{c346}, \cite{hoshi}.
It is known, that the above polynomial is irreducible for $m\neq -8,-3,0,5$
(see \cite{gras}), totally real and its cyclic Galois group is generated
by $x\rightarrow \frac{x-1}{x+2}$. Its cubic subfield $L_m$
is generated by a root of
\[
g_m(x)=x^3-mx^2-(m+3)x-1,
\]
these are the so called {\bf simplest cubic fields} of D.Shanks \cite{shanks}.
Its quadratic subfield is $H_m=\Q(\sqrt{q_m})$ with $q_m=m^2+3m+9$.
Let $\varphi_m$ be a root of the cubic polynomial $g_m(x)$
and $\psi_m=\sqrt{q_m}$ or $(1+\sqrt{q_m})/2$ according as 
$q_m\equiv 2,3 \ (\bmod \ 4)$ or $q_m\equiv 1 \ (\bmod \ 4)$.
Considering $K_m$ as a composite field of $L_m$ and $H_m$ we investigated 
power integral bases in the order 
$\O_m=\Z[1,\varphi_m,\varphi_m^2,\psi_m,\varphi_m\psi_m,\varphi_m^2\psi_m]$ 
of the ring of integers $\Z_{K_m}$ of $K_m$, and showed \cite{gop} that
$\O_m$ is not monogenic.

Our main purpose is now to determine an integral basis of $K_m$
in a parametric form and to consider monogenity in $\Z_{K_m}$.

\section{Integral bases of simplest sextic fields}

Our Theorem \ref{th1} below gives an integral basis of 
$K_m$ in a parametric form and implies that the structure of this integral 
basis is periodic in $m$ with period length 36.

For given $r$ with $1\leq r\leq 36$ consider those values $m=r+36k$ ($k\in\Z$)
for which $m\neq -8,-3,0,5$ and $q_m=m^2+3m+9$ is square free.
We omit the $r$ divisible by 3, because for these
values of $r$ the $q_m=m^2+3m+9$ is never square free.
As above, $\alpha_m$ denotes a root of (\ref{f}).

\begin{thm}
\label{th1} 
An integral basis $B_m$ of $K_m=\Q(\alpha_m)$ 
is obtained by substituting
$x=\alpha_m$ in the following formulas.
The discriminant of $K_m$ is given by $D_{K_m}$ in each case.
\end{thm}

\noindent
$r=1$, $m=r+36k$, $m^2+3m+9$ square free\\
$\displaystyle{
B_m=\left\{1,x,x^2,\frac{1+x+x^3}{2},\frac{4+x+3x^2+x^4}{6},
\frac{11+3x+13x^2+6x^3+2x^4+x^5}{18}\right\}
}$
\\
$D_{K_m}=(m^2+3m+9)^5$

\vspace{0.5cm}
\noindent
$r=2$, $m=r+36k$, $m^2+3m+9$ square free\\
$\displaystyle{
B_m=\left\{1,x,x^2,x^3,\frac{1+x^3+x^4}{3},
\frac{2+7x+6x^2+2x^3+x^5}{9}\right\}
}$
\\
$D_{K_m}=2^6\cdot (m^2+3m+9)^5$

\vspace{0.5cm}
\noindent
$r=4$, $m=r+36k$, $m^2+3m+9$ square free\\
$\displaystyle{
B_m=\left\{1,x,x^2,\frac{1+x^2+x^3}{2},\frac{1+x+3x^2+x^4}{6},
\frac{8+3x+x^2+3x^3+2x^4+x^5}{18}\right\}
}$
\\
$D_{K_m}=(m^2+3m+9)^5$

\vspace{0.5cm}
\noindent
$r=5$, $m=r+36k$, $m^2+3m+9$ square free\\
$\displaystyle{
B_m=\left\{1,x,x^2,\frac{1+x+x^3}{2},\frac{1+3x^2+x^3+x^4}{6},
\frac{8+10x+3x^2+5x^3+x^5}{18}\right\}
}$
\\
$D_{K_m}=(m^2+3m+9)^5$

\vspace{0.5cm}
\noindent
$r=7,34$, $m=r+36k$, $m^2+3m+9$ square free\\
$\displaystyle{
B_m=\left\{1,x,x^2,x^3,\frac{1+x+x^4}{3},
\frac{5+3x+7x^2+2x^4+x^5}{9}\right\}
}$
\\
$D_{K_m}=2^6\cdot (m^2+3m+9)^5$

\vspace{0.5cm}
\noindent
$r=8$, $m=r+36k$, $m^2+3m+9$ square free\\
$\displaystyle{
B_m=\left\{1,x,x^2,\frac{1+x^2+x^3}{2},\frac{4+3x+x^3+x^4}{6},
\frac{5+13x+8x^3+x^5}{18}\right\}
}$
\\
$D_{K_m}=(m^2+3m+9)^5$

\vspace{0.5cm}
\noindent
$r=10,19$, $m=r+36k$, $m^2+3m+9$ square free\\
$\displaystyle{
B_m=\left\{1,x,x^2,x^3,\frac{1+x+x^4}{3},
\frac{2+3x+4x^2+6x^3+2x^4+x^5}{9}\right\}
}$
\\
$D_{K_m}=2^6\cdot (m^2+3m+9)^5$

\vspace{0.5cm}
\noindent
$r=11$, $m=r+36k$, $m^2+3m+9$ square free\\
$\displaystyle{
B_m=\left\{1,x,x^2,x^3,\frac{1+x^3+x^4}{3},
\frac{2+7x+6x^2+2x^3+x^5}{9}\right\}
}$
\\
$D_{K_m}=2^6\cdot (m^2+3m+9)^5$

\vspace{0.5cm}
\noindent
$r=13$, $m=r+36k$, $m^2+3m+9$ square free\\
$\displaystyle{
B_m=\left\{1,x,x^2,\frac{1+x+x^3}{2},\frac{4+x+3x^2+x^4}{6},
\frac{8+12x+x^2+3x^3+2x^4+x^5}{18}\right\}
}$
\\
$D_{K_m}=(m^2+3m+9)^5$

\vspace{0.5cm}
\noindent
$r=14,23$, $m=r+36k$, $m^2+3m+9$ square free\\
$\displaystyle{
B_m=\left\{1,x,x^2,x^3,\frac{1+x^3+x^4}{3},
\frac{8+x+3x^2+5x^3+x^5}{9}\right\}
}$
\\
$D_{K_m}=2^6\cdot (m^2+3m+9)^5$

\vspace{0.5cm}
\noindent
$r=16$, $m=r+36k$, $m^2+3m+9$ square free\\
$\displaystyle{
B_m=\left\{1,x,x^2,\frac{1+x^2+x^3}{2},\frac{1+x+3x^2+x^4}{6},
\frac{5+3x+16x^2+2x^4+x^5}{18}\right\}
}$
\\
$D_{K_m}=(m^2+3m+9)^5$

\vspace{0.5cm}
\noindent
$r=17$, $m=r+36k$, $m^2+3m+9$ square free\\
$\displaystyle{
B_m=\left\{1,x,x^2,\frac{1+x+x^3}{2},\frac{1+3x^2+x^3+x^4}{6},
\frac{5+13x+9x^2+8x^3+x^5}{18}\right\}
}$
\\
$D_{K_m}=(m^2+3m+9)^5$

\vspace{0.5cm}
\noindent
$r=20$, $m=r+36k$, $m^2+3m+9$ square free\\
$\displaystyle{
B_m=\left\{1,x,x^2,\frac{1+x^2+x^3}{2},\frac{4+3x+x^3+x^4}{6},
\frac{11+7x+6x^2+2x^3+x^5}{18}\right\}
}$
\\
$D_{K_m}=(m^2+3m+9)^5$

\vspace{0.5cm}
\noindent
$r=22,31$, $m=r+36k$, $m^2+3m+9$ square free\\
$\displaystyle{
B_m=\left\{1,x,x^2,x^3,\frac{1+x+x^4}{3},
\frac{8+3x+x^2+3x^3+2x^4+x^5}{9}\right\}
}$
\\
$D_{K_m}=2^6\cdot (m^2+3m+9)^5$

\vspace{0.5cm}
\noindent
$r=25$, $m=r+36k$, $m^2+3m+9$ square free\\
$\displaystyle{
B_m=\left\{1,x,x^2,\frac{1+x+x^3}{2},\frac{4+x+3x^2+x^4}{6},
\frac{5+3x+7x^2+2x^4+x^5}{18}\right\}
}$
\\
$D_{K_m}=(m^2+3m+9)^5$

\vspace{0.5cm}
\noindent
$r=26,35$, $m=r+36k$, $m^2+3m+9$ square free\\
$\displaystyle{
B_m=\left\{1,x,x^2,x^3,\frac{1+x^3+x^4}{3},
\frac{5+4x+8x^3+x^5}{9}\right\}
}$
\\
$D_{K_m}=2^6\cdot (m^2+3m+9)^5$

\vspace{0.5cm}
\noindent
$r=28$, $m=r+36k$, $m^2+3m+9$ square free\\
$\displaystyle{
B_m=\left\{1,x,x^2,\frac{1+x^2+x^3}{2},\frac{1+x+3x^2+x^4}{6},
\frac{11+3x+4x^2+6x^3+2x^4+x^5}{18}\right\}
}$
\\
$D_{K_m}=(m^2+3m+9)^5$

\vspace{0.5cm}
\noindent
$r=29$, $m=r+36k$, $m^2+3m+9$ square free\\
$\displaystyle{
B_m=\left\{1,x,x^2,\frac{1+x+x^3}{2},\frac{1+3x^2+x^3+x^4}{6},
\frac{11+7x+15x^2+2x^3+x^5}{18}\right\}
}$
\\
$D_{K_m}=(m^2+3m+9)^5$

\vspace{0.5cm}
\noindent
$r=32$, $m=r+36k$, $m^2+3m+9$ square free\\
$\displaystyle{
B_m=\left\{1,x,x^2,\frac{1+x^2+x^3}{2},\frac{4+3x+x^3+x^4}{6},
\frac{8+x+3x^2+5x^3+x^5}{18}\right\}
}$
\\
$D_{K_m}=(m^2+3m+9)^5$.

\vspace{1cm}

\noindent
{\bf Proof} \\
\noindent
The discriminant of the polynomial (\ref{f}) is
\begin{equation}
D(f_m)=6^6\cdot q_m^5
\label{pp}
\end{equation}
with $q_m=m^2+3m+9$.
$\alpha_m$ denotes a root of $f_m$. 
Let ${\O_m}=\Z[\alpha_m]$. 
According to \cite{pz}, \cite{cohen} for a prime $p$
we call the order ${\O_m}$ 
$p$-maximal, if $p$ does not divide the index $[\Z_{K_m}^{+}:\O_m^{+}]$.

First we show that for all primes $p\neq 2,3$ dividing $q_m$, 
the order $\O_m$
is $p$-maximal. 
For this purpose we use Dedekind test, cf. \cite{pz}, \cite{cohen}.
For any polynomial $u(x)\in\Z[x]$, we set $\overline{u}(x)=u(x) (\bmod \; p\Z[x])$.

It is known (see \cite{gaal}), that
\[
f_m(x)\equiv \left( x- \frac{m}{3} \right)^6 \;\; (\bmod \; q_m\Z[x]).
\]
(Note that 3 does not divide $q_m$, therefore it is invertible mod $q_m$, as well as
mod $p$.)
We obtain 
\[
\overline{f_m}(x)=\left( x- \frac{m}{3} \right)^6.
\]
We set 
\[
g(x)=x- \frac{m}{3}
\]
and 
\[
h(x)=\frac{\overline{f_m}(x)}{\overline{g}(x)}=\left( x- \frac{m}{3} \right)^5.
\]
According to Dedekind test we construct
\[
t(x)=\frac{1}{p}(g(x)h(x)-f_m(x))=
\frac{1}{p}\left(\left(x- \frac{m}{3}\right)^6-f_m(x) \right).
\]
The order is $p$-maximal, if $\gcd(\overline{t}(x),\overline{g}(x),\overline{h}(x))=1$
in $\Z_p[x]$.
We have
\[
\left(x- \frac{m}{3}\right)^6-f_m(x)=
\frac{1}{3^6}(m^2+3m+9)\cdot \ell(x)
\]
with
\[
\ell(x)=1215x^4-540mx^3+1620x^3+135m^2x^2-405mx^2+54m^2x-18m^3x
\]
\[
-486x-3m^3+m^4+27m-81,
\]
satisfying
\[
\ell(x)=-3\cdot \left( x- \frac{m}{3} \right) \cdot 
(405x^3-45mx^2+540x^2+45mx+30m^2x+4m^3+33m^2-162)
\]
\[
+
(5(m+6)(m-3)(m^2+3m+9)+3^6).
\]
Indeed, by $p|(m^2+3m+9)$, $p\neq 2,3$ we obtain 
$\gcd(\overline{t}(x),\overline{g}(x),\overline{h}(x))=3^6$
which is a unit in $\Z_p[x]$.

\vspace{0.5cm}

We conclude, that only the primes 2 and 3 may occur in the denominators
of the basis elements. 

For all $r=1,\ldots,36$ such that $r^2+3r+9$ square free (24 possible values,
in case $r=5$ we consider $5+2\cdot 36=77$ instead of 5)
we may describe the types of an integral basis. We have to prove that for any
$m=r+36k$ this integral basis has the same structure like for $m=r$
(assuming both $r$ and $r+36k$ are square free).

\noindent
{\bf a)} Given $r$, using symmetrical polynomials we can prove that the basis elements
are algebraic integers for all $m=r+36k$. This is shown by explicitly calculating 
in all cases according to $r$ the minimal polynomial of all basis elements
and by checking that the coefficients are indeed integers.

\noindent
{\bf b)} If $m^2+3m+9$ is square free, then 3 does not divide it.
The square root of the discriminant (\ref{pp}) of $f_m$ is divisible by $3^3$.
For each type of our integral basis the product of the denominators 
is divisible by $3^3$, hence the discriminant of our integral basis is correct
according to the 3-factors.

\noindent
{\bf c)} 2 does not divide $m^2+3m+9$. 
The square root of the discriminant (\ref{pp}) of $f_m$ is divisible by $2^3$.
For certain types of the integral basis the product of the denominators 
is divisible by $2^3$. In these cases the discriminant of our integral basis 
is correct according to the 2-factors. 
For the other types of our integral bases the denominators are
not divisible by 2. In these cases (for these given values of $r$)
we test if some linear combination of the basis elements is divisible by 2.
If it were possible then also an element
\[
\beta=\frac{a_0+a_1\alpha+\ldots+a_5\alpha^5}{2}
\]
existed, which is algebraic integer. 
We set $m=r+4k$ and let $a_0,a_1,\ldots,a_5$ run through the residue
classes modulo 2, and $k$ run through the residue classes modulo $2^6$.
We calculate the minimal polynomial of $\beta$ (in the denominators of
some terms $2^6$ occur) and find that $\beta$ is never an algebraic integer
and this property only depends on the congruence behavior of $r$ 
modulo 4 and not on $k$. This was shown by explicitly calculating 
and checking the
minimal polynomial of $\beta$ for all given values of $r$ and for all
possible values of $a_0,\ldots,a_5,k$.  Hence also in these cases
the discriminant of our integral basis is correct
according to the 2-factors.

The discriminant $D_{K_m}$ is calculated from the discriminant of $f_m(x)$ and
the denominators of the elements of the integral basis.\\
$\Box$

\vspace{1cm}

\section{Monogenity of simplest sextic fields}

Using the integral basis we constructed
we managed to explicitly calculate the corresponding index form.
Note that it was a very complicated calculation because the index form is of degree 15 in 
the variables $x_2,\ldots,x_6$ with coefficients depending on $m$.
Investigating the factors of the index form we obtain the following statement.

\begin{thm}
\label{th2}
Let $m\neq -8,-3,0,5$ be an integer such that $q_m=m^2+3m+9$ is square free.
Let $K_m=\Q(\alpha_m)$ where $\alpha_m$ is a root of the polynomial (\ref{f}).
The field $K_m$ is not monogenic except for $m=-4,-2,-1,1$. 
For $m=-4,1$ and for $m=-2,-1$ the fields coincide.\\
For $m=1$ the field $K_m$ has discriminant $D_{K_m}=371293$ and
integral basis 
\[
B_1=\left\{1,x,x^2,\frac{1+x+x^3}{2},\frac{4+x+3x^2+x^4}{6},
\frac{11+3x+13x^2+6x^3+2x^4+x^5}{18}\right\}
\]
with $x=\alpha_1$ a root of (\ref{f}).
In this integral basis the coordinates $(y_1,y_2,y_3,y_4,y_5,y_6)$
of the generators of power integral bases
(up to sign) are given by\\
$(y_2,y_3,y_4,y_5,y_6)=
( 0 , 6 , 4 , 0 , -1),
( 1 , 1 , -2 , -2 , 1),\\
( 2 , -3 , -3 , -1 , 1),
( 4 , -2 , -6 , -3 , 2),
( 6 , -21 , -30 , -11 , 10),\\
( 7 , -10 , -20 , -9 , 7),
( 7 , -8 , -15 , -6 , 5),
( 8 , 1 , -10 , -7 , 4),\\
( 8 , 4 , -10 , -16 , 7),
( 8 , -3 , -9 , -4 , 3),
( 8 , -2 , -11 , -6 , 4),\\
( 9 , -5 , -14 , -7 , 5),
( 9 , -1 , -15 , -10 , 6),
( 9 , 3 , -16 , -13 , 7),\\
(  9, 1, -10, -7, 4),
(  9, 5, -11, -10, 5),
( 13 , -1 , -17 , -12 , 7),\\
( 15 , -5 , -20 , -10 , 7),
( 16 , -7 , -25 , -13 , 9),
( 17 , -4 , -24 , -14 , 9),\\
( 24 , -6 , -35 , -20 , 13),
( 24 , -12 , -39 , -20 , 14),
(  24 , -10 , -34 , -17 , 12),\\
( 28 , -102 , -147 , -54 , 49),
( 32 , -15 , -48 , -24 , 17),
( 33 , -5 , -45 , -27 , 17),\\
( 41 , -8 , -54 , -31 , 20),
( 41 , -14 , -58 , -31 , 21),
( 49 , -23 , -71 , -35 , 25),\\
( 57 , -9 , -75 , -44 , 28),
( 59 , -2 , -86 , -57 , 34),
( 105 , -54 , -157 , -76 , 55),\\
( 146 , -33 , -198 , -112 , 73),
( 246 , 13 , -359 , -250 , 145),\\
( 517 , -268 , -774 , -374 , 271),
(723 , -155 , -970 , -551 , 358),
$\\
independently of $y_1$.

\vspace{0.5cm}

For $m=-1$ the field $K_{-1}$ has discriminant $D_{K_{-1}}=1075648$ and
integral basis 
\[
B_{-1}=\left\{1,x,x^2,x^3,\frac{1+x^3+x^4}{3},
\frac{5+4x+8x^3+x^5}{9}\right\}
\]
with $x=\alpha_{-1}$ a root of (\ref{f}).
In this integral basis the coordinates $(y_1,y_2,y_3,y_4,y_5,y_6)$ 
of the generators of power integral bases
(up to sign) are given by\\
$(y_2,y_3,y_4,y_5,y_6)=
( 0, 0, 2, 0, -1),
( 1 , 13 , 18 , -4 , -8),\\
( 1 , 9 , 22 , -3 , -10),
( 2 , 6 , 11 , -2 , -5),
( 2 , 73 , 119 , -23 , -53),\\
( 2 , 9 , 11 , -3 , -5),
( 3 , 35 , 56 , -11 , -25),
( 3 , 3 , 2 , -1 , -1),\\
( 3 , 4 , 4 , -1 , -2),
( 3 , 11 , 37 , -4 , -17),
(5 , 10 , 11 , -3 , -5),\\
( 13 , 23 , 22 , -7 , -10),
( 27 , 96 , 112 , -29 , -50),
( 27 , 39 , 37 , -12 , -17),\\
$
independently of $y_1$.
\end{thm}

\vspace{1cm}

\noindent
{\bf Proof} \\
\noindent
Let $\alpha_m=\alpha_m^{(1)}$ be a root of (\ref{f}) and for $i=2,\ldots,6$
set 
\[
\alpha_m^{(i)}=\frac{\alpha_m^{(i-1)}-1}{\alpha_m^{(i-1)}+2}.
\]
Denote by $\{b_1=1,b_2,\ldots,b_6\}$ the integral basis of $K_m$ 
constructed in Theorem \ref{th1} and let
$b_i^{(j)}$ be the conjugate of $b_i$ corresponding to $\alpha_m^{(j)}$.
Let 
\[
L^{(j)}(\underline{X})=X_1+X_2b_2^{(j)}+\ldots +X_6b_6^{(j)}
\]
for $j=1,\ldots,6$. The index form corresponding to the integral basis
has three factors in our case:
\[
F_1(\underline{X})=
(L^{(1)}(\underline{X})-L^{(2)}(\underline{X}))
(L^{(2)}(\underline{X})-L^{(3)}(\underline{X}))
(L^{(3)}(\underline{X})-L^{(4)}(\underline{X}))
\]
\[
\cdot (L^{(4)}(\underline{X})-L^{(5)}(\underline{X}))
(L^{(5)}(\underline{X})-L^{(6)}(\underline{X}))
(L^{(6)}(\underline{X})-L^{(1)}(\underline{X}))
\]
\[
=N_{K/\Q}(L^{(1)}(\underline{X})-L^{(2)}(\underline{X})),
\]

\[
F_2(\underline{X})=
(L^{(1)}(\underline{X})-L^{(3)}(\underline{X}))
(L^{(2)}(\underline{X})-L^{(4)}(\underline{X}))
(L^{(3)}(\underline{X})-L^{(5)}(\underline{X}))
\]
\[
\cdot (L^{(4)}(\underline{X})-L^{(6)}(\underline{X}))
(L^{(5)}(\underline{X})-L^{(1)}(\underline{X}))
(L^{(6)}(\underline{X})-L^{(2)}(\underline{X}))
\]
\[
=N_{K/\Q}(L^{(1)}(\underline{X})-L^{(3)}(\underline{X})).
\]
and
\[
F_3(\underline{X})=
(L^{(1)}(\underline{X})-L^{(4)}(\underline{X}))
(L^{(2)}(\underline{X})-L^{(5)}(\underline{X}))
(L^{(3)}(\underline{X})-L^{(6)}(\underline{X}))
\]
\[
=N_{K/\Q}(L^{(1)}(\underline{X})-L^{(4)}(\underline{X})).
\]
These are polynomials
in $(X_2,\ldots,X_6)$, not depending on $X_1$.

In each case of our integral basis we 
explicitly calculated these polynomials.
This was an extensive calculation performed by Maple.
In all the 24 possible values of $r$ we calculated
explicitly $F_1(\underline{X}), F_2(\underline{X}), F_3(\underline{X})$.
These forms are of degree 6,6,3, respectively, in 5 variables,
having about 1500 coefficients depending on $m$.
Therefore their explicit form cannot be listed or considered
without a computer algebra system.

We have $D_{K_m}=2^{2\ell} (m^2+3m+9)^5$ with $\ell=0,3$.
Using a very careful calculation we could also factorize 
these polynomials by Maple and found that 
\[
F_i(\underline{X})=C_i \cdot G_i(\underline{X})\;\;\; (i=1,2,3),
\]
with
\[
C_1=C_2=m^2+3m+9,\;\;\; C_3=2^{\ell} \sqrt{m^2+3m+9},
\]
where the polynomials $G_i(\underline{X})$ turned out to have 
integer coefficients (here we used the congruence behavior of $m$ modulo 36).
By 
\[
\sqrt{D_{K_m}}=C_1\cdot C_2\cdot C_3
\]
the index form equation corresponding to the given basis 
is equivalent to
\[
G_1(\underline{x})\cdot G_2(\underline{x})\cdot G_3(\underline{x})=\pm 1,
\;\;  x_2,\ldots,x_6\in\Z.
\]
Therefore the existence of a power integral basis implies that 
there exist $x_2,\ldots,x_6\in\Z$ with
\[
G_i(\underline{x})=\pm 1 \;\; {\rm for}\;\; i=1,2,3.
\]
We observe that 
\[
q_m=m^2+3m+9\; | \; 27G_1(\underline{X})+G_2(\underline{X}).
\]
Hence for any solution $x_2,\ldots,x_6\in\Z$ of the index form equation
we have
\[
q_m=m^2+3m+9 \; |\; \pm 27\pm 1.
\]
This is only possible for $m=-4,-2,-1,1$, therefore
$K_m$ is not monogenic for $m\not=-4,-2,-1,1$.

Considering the fields $K_m$ for $m=-4,1$ we obtain that they have
the same quadratic and cubic subfields, hence they coincide.
We obtain the same for $m=-2,-1$.
(In general, we know from \cite{hoshi} that 
$K_m$ and $K_{-m-3}$ are the same algebraic number fields.)

The totally real cyclic sextic fields with discriminants
$371293$ and $1075648$ were investigated in \cite{gaalc6}.
These fields have the same subfields as $K_1$ and $K_{-1}$,
respectively, hence they are the same fields. 
In \cite{gaalc6}
all generators of power integral bases of these fields were
calculated in a different integral basis. 
The $(y_2,\ldots,y_6)$ in our Theorem are obtained by 
converting the coordinates of those elements into our 
integral bases.\\
$\Box$

\vspace{1cm}

\section{Computational remarks}

All calculations mentioned in the paper were implemented
in Maple \cite{maple} and executed on an average laptop.
In the proof of Theorem \ref{th1} the test to check if
the discriminant of our integral basis has the correct 2 factor
had $2^6\cdot 2^6=4096$ steps which took just some seconds in each case.
The explicit calculation of the factors of the index form in 
the proof of Theorem \ref{th2} also only took some seconds,
involving factorization, in each case of our integral basis.
All together the CPU time for the proofs was very short but 
we had to organize the calculations very efficiently
because the size of the polynomials were on the limits of the capacity of the
Maple system.
For example the polynomials $F_i$ in the proof of Theorem \ref{th2}
had 1500-1800 terms, hence the index form had about $4\cdot 10^9$ terms.

\end{document}